\documentclass[12pt]{amsart}
\usepackage{latexsym}
\usepackage{a4}
\usepackage{amssymb}
\usepackage{amsmath}
\usepackage[english]{babel}
\usepackage{tikz}
\usepackage{bm}

\newtheorem{thm}{Theorem}[section]
\newtheorem{pro}[thm]{Proposition}
\newtheorem{lem}[thm]{Lemma}

\newtheorem{cor}[thm]{Corollary}
\theoremstyle{definition}
\newtheorem{exa}[thm]{Example}
\newtheorem{de}[thm]{Definition}
\newtheorem{rem}[thm]{Remark}
\numberwithin{equation}{section}

\newcommand{\N}{\mathbb{N}}
\newcommand{\Lg}{\mathcal{L}}
\newcommand{\Rg}{\mathcal{R}}

\newcommand{\Hg}{\mathcal{H}}

\newcommand{\Con}{\operatorname{Con}}

\newcommand{\Pol}{\operatorname{Pol}}

\newcommand{\Term}{\operatorname{Term}}
\newcommand{\ar}{\operatorname{ar}}
\newcommand{\bS}{\mathbf{S}}
\newcommand{\bA}{\mathbf{A}}
\newcommand{\bG}{\mathbf{G}}
\newcommand{\oa}{\mathbf{a}}
\newcommand{\ob}{\mathbf{b}}
\newcommand{\oc}{\mathbf{c}}
\newcommand{\od}{\mathbf{d}}

\newcommand{\oy}{\mathbf{y}}
\newcommand{\oz}{\mathbf{z}}

\title{Commutators in completely simple semigroups}
\author{Jelena Radovi\'{c}}
\author{Neboj\v{s}a Mudrinski}

\thanks{Supported by the Ministry of Science, Education and Technological Development of the Republic of Serbia (Grant No. 451-03-68/2022-14/200125),
and the Ministry of Scientific and Technological Development, Higher Education and Information Society of the Republic of Srpska (Grant No. 19.041/68-13-27/2023.) }
\date{\today}

\begin{document}
\bibliographystyle{amsalpha}

\maketitle

\begin{center}
\begin{small}
Jelena Radovi\'{c}\footnote{Corresponding author}\\
Department of Mathematics\\
University of East Sarajevo\\
71123 East Sarajevo\\
Bosnia and Herzegovina\\
ORCID: 0000-0002-1023-4463\\
{\tt jelena.radovic@ff.ues.rs.ba}\\[5mm]
\end{small}
\end{center}
\begin{center}
\begin{small}
Neboj\v{s}a Mudrinski\\
Department of Mathematics and Informatics\\
Faculty of Sciences\\
University of Novi Sad\\
21000 Novi Sad\\
Serbia\\
ORCID: 0000-0001-9830-6603\\
{\tt nmudrinski@dmi.uns.ac.rs}\\[5mm]
\end{small}
\end{center}

\begin{abstract}
We obtain a characterization of the binary commutator on completely simple semigroups, using their Rees matrix representation. Consequently, we prove that a regular semigroup is nilpotent (solvable) if and only if it is simple, and all its $\mathcal{H}$-classes are nilpotent (solvable) groups.
\\
\bigskip

\end{abstract}

\section{Motivation and Results}

Let $\bS=(S,\cdot)$ denote a semigroup, that is, an algebra with an associative binary operation $\cdot$ on a nonempty set $S$. In general, we say that an algebra $\bA$ is \emph{abelian} (TC algebra), if for every $(n+1)$-ary term $t=t(x,\oy)$, and for every $a,b\in A$, $\oc,\od\in A^n$ we have $t(a,\oc)=t(a,\od)\; \implies \; t(b,\oc)=t(b,\od)$. This has been introduced by R. McKenzie in \cite{MK:OMLFV}, and it can be translated to the equality $[1_A,1_A]=0_A$ (see Definition \ref{McMcDefn4.150}), where $1_A$ denotes the full relation and $0_A$ denotes the equality on the set $A$. R. McKenzie \cite{MK:TNoNIM} has described abelian semigroups of finite exponent ($x^n\approx x^{2n}$). W. Taylor \cite{Taylor} has also studied abelian semigroups. R.J. Warne has obtained a complete characterization of abelian semigroups, and more specifically, regular abelian  semigroups, in \cite{Warne,Warne2}. A semigroup $\bS$ is \emph{regular} if for every $a\in S$ there exists an $x\in S$ such that $axa=a$. In this paper we study nilpotency ($[1_S,\dotsc,[1_S,1_S]]=0_S$) and solvability (sequence $[1_S,1_S]$, $[[1_S,1_S],[1_S,1_S]]$, $\dotsc$ collapses to $0_S$) of regular semigroups. However, due to \cite[Proposition 3.5]{RM23}, this question reduces to investigation of nilpotency and solvability in completely simple semigroups. For this purpose we provide a characterization of the binary commutator in completely simple semigroups using their Rees matrix representation, and corresponding representation of congruences by the linked triples. Throughout this paper we will exclusively use the so called normalized Rees matrix representations $\mathcal{M}[G;I,\Lambda;P]$, as described in Proposition \ref{HowieTe3.3.1}, without explicit referencing. Due to the isomorphism between the lattice of normal subgroups of the group $\bG$ and its congruence lattice, in place of $N_{\rho}$ in the standard linked triple $(\rho_I,N_{\rho},\rho_{\Lambda})$ (see \cite[page 91]{Howie}), we can consider the group congruence $\rho_G$ as a part of the linked triple. Here, the congruence $\rho\in\Con(\bS)$ determines the group congruence $\rho_G$ on $\bG$ by $a\;\rho_G\;b \iff (1,a,1)\;\rho\;(1,b,1)$. The main result of this note is the following statement.

\begin{thm}\label{ThmComm}
Let $\bS = \mathcal{M}[G;I,\Lambda;P]$ be a completely simple semigroup, and let $\rho,\sigma$ be congruences on $\bS$. The commutator $[\rho,\sigma]$ corresponds to the linked triple $(0_I,[\rho_G,\sigma_G]\vee \Theta_{\rho,\sigma},0_{\Lambda})$.
\end{thm}

Here $\Theta_{\rho,\sigma}$ denotes the congruence on $\bG$ determined by the matrix $P$ and components of $\rho$ and $\sigma$ on the index sets $I$ and $\Lambda$, see Definition \ref{thetarhosigma}. Note that if we consider an expanded completely simple semigroup $(S,\cdot,^{-1})$ by the unary operation $^{-1}$ that maps an element to its inverse in the corresponding maximal subgroup, the commutator will be the same as the commutator in $(S,\cdot)$, by Proposition \ref{PropAFu} because $^{-1}$ is compatible with all congruences of $(S,\cdot)$ (\cite[Exercise 4.7.2]{Howie}). As a corollary of Theorem \ref{ThmComm}, we deduce that the binary commutator on completely simple semigroups behaves as in congruence modular algebras: every congruence above it centralizes the congruences in the commutator (Proposition \ref{Cor1}), the commutator is symmetric (Proposition \ref{Cor[a,b]=[b,a]}) and agrees with the factor congruences (Proposition \ref{PropFactor}) and the join of congruences (Proposition \ref{Cor[a,bigveeb]=bigvee[a,b]}). The stated description of the commutator enables us to prove in Proposition \ref{PropNilpSolv} the equivalence between the nilpotency (solvability) of a completely simple semigroup $\bS=\mathcal{M}[G;I,\Lambda;P]$ and the group $\bG$, that is isomorphic to its maximal subgroup. This leads to the characterization of nilpotency and solvability in a regular semigroup in the following theorem.

\begin{thm}\label{ThmNilpSolv}
Let $\bS$ be a regular semigroup. Then:
\begin{enumerate}
\item[(a)] $\bS$ is nilpotent if and only if it is simple and all its $\mathcal{H}$-classes are nilpotent groups;
\item[(b)] $\bS$ is solvable if and only if it is simple and all its $\mathcal{H}$-classes are solvable groups.
\end{enumerate}
\end{thm}

\section{Preliminaries From Semigroup Theory}

We will use the standard notation $\Term(\bS)$ and $\Pol(\bS)$ for the sets of all term functions and polynomials, respectively, over a semigroup $\bS$. In particular, by $\Term_k(\bS)$ and $\Pol_k(\bS)$, respectively, we will denote the sets of all term functions and polynomials of arity $k$, $k\in\N$. Recall also that a term $t=t(x_1,\dotsc,x_n)$ of semigroup type is a word over the alphabet $\{x_1,\dotsc,x_n\}$, while a polynomial term $p=p(x_1,\dotsc,x_n)$ of semigroup type is a word over the alphabet $\{x_1,\dotsc,x_n,c_1,\dotsc,c_m\}$, for some constants $c_1,\dotsc,c_m$. We will also use the usual notation $\Con(\bS)$ for the lattice of all congruences on a semigroup $\bS$.

A semigroup $\bS$ is \textit{simple} if it has no proper ideals (see \cite[page 5]{Howie}). \emph{Idempotents} in semigroup $\bS$ are elements of $S$ that satisfy the identity $x^2=x$. A \emph{primitive idempotent} is an element $e$ that uniquely satisfies the equations $ex=xe=x=x^2$ on $x$. A simple semigroup that contains primitive idempotents is called a \emph{completely simple semigroup}.

\begin{pro}\textup{(\cite{Howie} Theorem 3.3.1, Theorem 3.4.2)}\label{HowieTe3.3.1}
Let $\bG$ be a group with the identity element $e$, let $I$ and $\Lambda$ be non-empty sets such that there exists an element $1\in I\cap\Lambda$, and let $P=[p_{\lambda i}]$ be a $\Lambda\times I$ matrix with entries in $G$ that satisfies the property $p_{1 i}=e=p_{\lambda 1}$ for every $i\in I,\lambda\in\Lambda$. Let $S=I\times G\times\Lambda$, and define multiplication on $S$ by
\begin{equation*}
(i,g,\lambda) \cdot (j,h,\mu) = (i,gp_{\lambda j} h,\mu).
\end{equation*}
Then $\mathcal{M}[G;I,\Lambda;P]=(S,\cdot)$ is a completely simple semigroup. Conversely, every completely simple semigroup is isomorphic to a semigroup constructed in this way.
\end{pro}

Representation of a completely simple semigroup by Rees matrix semigroups gives a description of its congruences in terms of linked triples. Every congruence $\rho$ on a completely simple semigroup $\bS=\mathcal{M}[G;I,\Lambda;P]$ uniquely determines a \emph{linked triple} $(\rho_I,N_{\rho},\rho_{\Lambda})$, where $\rho_I$ and $\rho_{\Lambda}$ are equivalences on $I$ and $\Lambda$, respectively, defined by
\begin{align*}
\rho_I &= \{(i,j)\in I^2 : (\forall \lambda\in\Lambda) (i,p_{\lambda i}^{-1},\lambda)\; \rho \; (j,p_{\lambda j}^{-1},\lambda)\},\\
\rho_{\Lambda} &= \{(\lambda,\mu)\in \Lambda^2 : (\forall i\in I) (i,p_{\lambda i}^{-1},\lambda) \;\rho \; (i,p_{\mu i}^{-1},\mu)\},
\end{align*}
while $N_{\rho} = \{g\in G: (1,g,1)\; \rho\; (1,e,1)\}$ is a normal subgroup of $\bG$. The correspondence $:\rho\rightarrow (\rho_I,N_{\rho},\rho_{\Lambda})$ is an order preserving bijection (see \cite[Theorem 3.5.9]{Howie}). Since we use linked triples that have changed component $N_{\rho}$ to $\rho_G$, this statement can be reformulated as in the following proposition.

\begin{pro}\label{PropCongOrd}
Let $\bS=\mathcal{M}[G;I,\Lambda;P]$ be a completely simple semigroup. The correspondence $:\rho\mapsto (\rho_I,\rho_G,\rho_{\Lambda})$ is an order preserving bijection between the lattice $\Con(\bS)$ and the lattice of all linked triples.
\end{pro}

The following proposition can be found in \cite[Proposition 2.10]{RM23}. We will use it without explicit referencing.

\begin{pro}\label{PropCongGS}
Let $\bS=\mathcal{M}[G;I,\Lambda;P]$ be a completely simple semigroup and let $\rho\in\Con(\bS)$. Then for every $(i,g,\lambda),(j,h,\mu)\in S$ we have $(i,g,\lambda)\;\rho\;(j,h,\mu)$ if and only if $i\;\rho_I\;j$, $g \;\rho_G\; h$ and $\lambda\;\rho_{\Lambda}\;\mu$.
\end{pro}

Here we recall some basic notions about regular semigroups (see \cite{CP:ATS,Howie}). Green's equivalences $\mathcal{R}$, $\mathcal{L}$ and $\mathcal{H}$ on an arbitrary regular semigroup are not necessarily congruences.  However, if $\bS=\mathcal{M}[G;I,\Lambda;P]$ is a completely simple semigroup, then $\Lg$ and $\Rg$ are congruences (\cite[Exercise 2.7.9]{CP:ATS}), and consequently, their intersection $\Hg$ is also a congruence. In the language of linked triples, congruence $\mathcal{H}$ is uniquely determined with the linked triple $(0_I,1_G,0_{\Lambda})$, while congruences $\Lg$ and $\Rg$ are determined by linked triples $(1_I,1_G,0_{\Lambda})$ and $(0_I,1_G,1_{\Lambda})$, respectively. Moreover, every $\Hg$-class is a maximal subgroup of $\bS$, isomorphic to the group $\bG$. A semigroup with a unary operation $^{-1}:a\mapsto a^{-1}$ that satisfies $(a^{-1})^{-1}=a$, $aa^{-1}a=a$ and $aa^{-1}=a^{-1}a$ is called a \emph{completely regular semigroup}.

\section{Centralizing Condition}

For $n\in\N$ and $A\neq\emptyset$ we denote an arbitrary element of $A^n$ by $\mathbf x$ and its coordinates by $x_1,\dots,x_n$, that is, ${\mathbf x}=(x_1,\dots,x_n)$. If $\bA$ is an algebra and $\theta$ a congruence of $\bA$ we write $\oa\;\theta\;\ob$ for $\oa,\ob\in A^n$, $n\in\N$ if $a_i\;\theta\; b_i$ for all $i\in\{1,\dots,n\}$.

\begin{de}\textup{(\cite{McMc}, Definition 4.148)}\label{McMcDefn4.148}
Let $\bA$ be an algebra, and $\alpha,\beta,\delta$ congruences on $\bA$. We say that \textit{$\alpha$ centralizes $\beta$ modulo $\delta$}, in abbreviation $C(\alpha,\beta;\delta)$, if for all $t\in \Term_{n+1}(\bA)$, $n\in\N$ and for all $a,b\in A$, $\oc,\od\in A^n$ such that $a\;\alpha\;b$, $\oc\;\beta\;\od$, the following implication is true:
\begin{equation}
t(a,\oc)\; \delta\; t(a,\od) \implies t(b,\oc)\;\delta\; t(b,\od). \label{TCalg}
\end{equation}
\end{de}

\begin{de}\textup{(\cite{McMc}, Definition 4.150)}\label{McMcDefn4.150}
Let $\bA$ be an algebra. \textit{The commutator $[\alpha,\beta]$} of congruences $\alpha,\beta\in\Con(\bA)$ is the least congruence $\delta$ such that $\alpha$ centralizes $\beta$ modulo $\delta$.
\end{de}

\begin{pro}\textup{(\cite{FM:CTCMV}, Proposition 3.4)} \label{CommutatorProperties}
Let $\bA$ be an algebra and let $\alpha,\beta,\alpha_1$ and $\beta_1$ be congruences on $\bA$. If $\alpha_1\leq\alpha$ and $\beta_1\leq\beta$, then $[\alpha_1,\beta_1]\leq[\alpha,\beta]$.
\end{pro}

\begin{pro}\textup{(cf. \cite{FM:CTCMV}, Proposition 4.2)} \label{PropCongModVar}
Let $\mathbf{A}$ be an algebra from a congruence modular variety and let $\alpha,\beta,\delta \in\Con(\mathbf{A})$. Then $C(\alpha,\beta;\delta)$ if and only if $[\alpha,\beta]\leq \delta$.
\end{pro}

\begin{pro}\label{PropAFu}
Let $\bA=(A,\mathcal{F})$ and $\mathcal{A}^u=(A,\mathcal{F}\cup\mathcal{U})$ be algebras with the same universe $A$, where $\mathcal{U} = \{f_i:i \in I\}$ is a set of operations $f_i$ of arity at most $1$. If $\Con(\bA)=\Con(\bA^u)$, then for every $\alpha,\beta,\delta\in\Con(\bA)$ we have $C(\alpha,\beta;\delta)_A$ if and only if $C(\alpha,\beta;\delta)_{A^u}$.
\end{pro}

\begin{proof}
Without loss of generality, we will give the proof for the case $\mathcal{U}=\{f\}$, $\ar(f)=1$. Let $\alpha,\beta,\delta$ be arbitrary congruences from $\Con(\bA)=\Con(\bA^u)$. From Definition \ref{McMcDefn4.148} one can easily see that $C(\alpha,\beta;\delta)_{A^u}$ implies $C(\alpha,\beta;\delta)_A$, since $\Term(\bA)\subseteq\Term(\bA^u)$. On the other hand, let $\delta\in\Con(\bA)$ be such that $C(\alpha,\beta;\delta)_A$. Let $t=t(x,\oy)$ be a term in $\Term_{n+1}(\bA^u)$, $n\in\N$ and let $a,b\in A$, $\oc,\od\in A^n$ be such that $a\;\alpha\;b$, $\oc\;\beta\;\od$, and $t(a,\oc)\;\delta\;t(a,\od)$. Define the term $s\in\Term_{2+2n}(\bA)$, as $s(x,w,y_1,\dotsc,y_n,z_1,\dotsc,z_n)$ in the following way: term $s$ is obtained by replacing every occurrence of $f(x)$ in the term $t$ by the new variable $w$, and also by replacing every occurrence of $f(y_i)$ in $t$ by the new variable $z_i$, for $i=1,\dotsc,n$. Then we have $t(x,\oy)=s(x,f(x),\oy,f(\oy))$ for every $x\in A$, $\oy\in A^n$, where we have denoted $f(\oy)=(f(y_1),\dotsc,f(y_n))$. Since $\Con(\bA)=\Con(\bA^u)$, relations $a\;\alpha\;b$ and $\oc\;\beta\;\od$ imply $f(a)\;\alpha\;f(b)$ and $f(\oc)\;\beta\;f(\od)$, respectively. Hence, we have $(f(a),\oc,f(\oc))\;\beta\;(f(a),\od,f(\od))$ and $(b,\oc,f(\oc))\;\beta\;(b,\od,f(\od))$. Therefore, using $C(\alpha,\beta;\delta)_A$ from $s(a,f(a),\oc,f(\oc))=t(a,\oc)\;\delta\;t(a,\od)=s(a,f(a),\od,f(\od))$ we obtain $s(b,f(a),\oc,f(\oc))\;\delta\;s(b,f(a),\od,f(\od))$. Now, define $r\in\Term_{2+2n}(\bA)$ by $r(w,x,\oy,\oz):=s(x,w,\oy,\oz)$. Again using the condition $C(\alpha,\beta;\delta)_A$, and $r(f(a),b,\oc,f(\oc))=s(b,f(a),\oc,f(\oc))\;\delta\;s(b,f(a),\od,f(\oy))=r(f(a),b,\od,f(\od))$, we obtain

$s(b,f(b),\oc,f(\oc))=r(f(b),b,\oc,f(\oc))\;\delta\;r(f(b),b,\od,f(\od))=s(b,f(b),\od,f(\od))$. Hence, we have $t(b,\oc)\;\delta\;t(b,\od)$, which completes our proof of $C(\alpha,\beta;\delta)_{A^u}$.
\end{proof}

Note that if $\mathcal{U}$ contains just constants, then $\Term(\bA^u)=\Pol(\bA)$, and we obtain the following lemma, stated in \cite[Excercises 4.156.2]{McMc} and proved in \cite{EA:TPFoCA}.

\begin{lem}\label{McMcExer2}
Let $\mathbf{A}$ be an algebra, and let
$\alpha,\beta,\delta\in\Con(\bA)$. Then $C(\alpha,\beta;\delta)$ if and only if for every polynomial $p\in\Pol_{n+1}(\bA)$, $n\in\N$ and for all $a,b\in A$, $\oc,\od\in A^n$, $n\in\N$ such that $a\;\alpha\; b$ and $\oc\;\beta\;\od$, the following implication is true:
\begin{equation}
p(a,\oc) \;\delta\; p(a,\od) \implies p(b,\oc)\;\delta\; p(b,\od). \label{PC}
\end{equation}
\end{lem}

In order to simplify checking the centralizing condition for semigroups, we generalize (TC1), (TC2) and (TC3) from \cite{Warne}. We will say that congruences $\alpha,\beta,\delta$ on a semigroup $\bS$ satisfy
\begin{enumerate}
\item[\textup{(C1)}] if $ac\,\delta\,ad \Rightarrow bc\,\delta\,bd$;
\item[\textup{(C2)}] if $ca\,\delta\,da \Rightarrow cb\,\delta\,db$;
\item[\textup{(C3)}] if $c_1ac_2\,\delta\,d_1ad_2 \Rightarrow c_1bc_2\,\delta\,d_1bd_2$;
\end{enumerate}
for all $a,b,c,d,c_1,d_1,c_2,d_2\in S$ such that $a\;\alpha\;b$, $c\;\beta\;d$, $c_1\;\beta\;d_1$, $c_2\;\beta\;d_2$.

\begin{pro}\label{TCOpsti}
Let $\bS$ be a semigroup, and let $\alpha,\beta,\delta\in\Con(\bS)$. Then we have $C(\alpha,\beta;\delta)$ if and only if $\alpha,\beta$ and $\delta$ satisfy {\rm (C1)}, {\rm (C2)} and {\rm (C3)}.
\end{pro}

\begin{proof}
$(\rightarrow)$ Assume that we have $C(\alpha,\beta;\delta)$. We obtain (C1), (C2) and (C3) using Definition \ref{McMcDefn4.148} for terms $t_1(x,y)=xy$, $t_2(x,y)=yx$ and $t_3(x,y,z)=yxz$, respectively. $(\leftarrow)$ Now assume that the conditions (C1), (C2) and (C3) are satisfied in the semigroup $\bS$. We proceed using Lemma \ref{McMcExer2}. Let $p\in\Pol_{n+1}(\bS)$, $n\in\N$ and let $a,b\in S$ and $\oc,\od\in S^n$ be such that $a\;\alpha\;b$ and $\oc\;\beta\;\od$. We note that for any polynomials $f\in\Pol_{n+1}(\bS)$ and $g\in\Pol_n(\bS)$, we
have $f(b,\oc)\;\beta\;f(b,\od)$ and $g(\oc)\;\beta\;g(\od)$. Assume that $p(a,\oc)\;\delta\; p(a,\od)$. We will prove that $p(b,\oc)\;\delta\; p(b,\od)$ by induction on the number $k$, $k\geq 0$ of appearances of $x$ in the corresponding word of the polynomial $p(x,\mathbf{y})$. First, let $k=0$, then $p(x,\oy)=q(\oy)$ where $q\in\Pol_n(\bS)$. Therefore, we have $q(\oc)\;\delta\;q(\od)$, which trivially implies $p(b,\oc)\;\delta\;p(b,\od)$. Now assume
that the statement is true for every polynomial $p\in\Pol_{n+1}(\bS)$, $n\in\N$, with less than $k$ occurrences of
$x$ in the corresponding polynomial word, $k\geq 1$ and let $p(x,\oy)$ be a polynomial over $\bS$, with exactly $k$ occurrences of $x$ in the corresponding word of $p(x,\oy)$. Let us observe the first occurrence of $x$ in the corresponding word of $p(x,\mathbf{y})$. We have the following three possibilities:

\noindent (i) $p(x,\mathbf{y})=x\cdot q(x,\mathbf{y})$, where $q$ has exactly $k-1$ occurrences of $x$ in the corresponding word of $q(x,\oy)$. Define the polynomial $r(x,\oy)=a\cdot q(x,\oy)$. Then we have $r(a,\oc) = a\cdot q(a,\oc) = p(a,\oc) \;\delta\; p(a,\od) = a\cdot q(a,\od) = r(a,\od)$. Since $r(x,\oy)$ has $k-1$ occurrences of $x$ in the corresponding word, we obtain $a\cdot q(b,\oc) = r(b,\oc) \;\delta\; r(b,\od) = a\cdot q(b,\od)$ by the inductive hypothesis. Hence we have $p(b,\oc) = b\cdot q(b,\oc) \;\delta\; b\cdot q(b,\od) = p(b,\od)$ using $q(b,\oc)\;\beta\; q(b,\od)$ and the condition (C1).

\noindent (ii) $p(x,\mathbf{y})=q(\mathbf{y})\cdot x$. Let us note that in this case $q$ has exactly $0$ occurrences of $x$ and our assumption simplifies to $q(\oc)\cdot a \;\delta\; q(\od)\cdot a$. Therefore, we have $p(b,\oc) = q(\oc)\cdot b \;\delta\; q(\od)\cdot b = p(b,\od)$ using $q(\oc)\;\beta\;q(\od)$ and the condition (C2).

\noindent (iii) $p(x,\oy)=q_1(\oy)\cdot x \cdot q_2(x,\oy)$, where $x$ does not occur in the corresponding word of $q_1$ and $q_2$ has exactly $k-1$ occurrences of $x$ in its corresponding word. Define the polynomial $s(x,\oy)=q_1(\oy)\cdot a\cdot q_2(x,\oy)$. Then we have
\begin{equation*}
s(a,\oc) = q_1(\oc)\cdot a\cdot q_2(a,\oc) = p(a,\oc) \;\delta\; p(a,\oc) = q_1(\od) \cdot a \cdot q_2(a,\od)=s(a,\od).
\end{equation*}
Since $s$ has $k-1$ occurrences of $x$ in the corresponding word, we obtain $q_1(\oc)\cdot a\cdot q_2(b,\oc) = s(b,\oc) \;\delta\; s(b,\od) = q_1(\od) \cdot a\cdot q_2(b,\od)$ by the inductive hypothesis. Hence we have $p(b,\oc) =
q_1(\oc) \cdot b \cdot q_2(b,\oc) \;\delta\; q_1(\od) \cdot b\cdot q_2(b,\od) = p(b,\od)$ using $q_1(\oc)\;\beta\;q_1(\od)$, $q_2(b,\oc)\;\beta\;q_2(b,\od)$ and the condition (C3).
\end{proof}

In the case of completely simple semigroups the condition $C(\rho,\sigma;\theta)$ can be simplified as follows.

\begin{pro}\label{TCuslovRMSC3}
Let $\bS$ be a completely simple semigroup and let $\rho,\sigma,\theta$ be congruences on $\bS$. Then the condition $C(\rho,\sigma;\theta)$ is true if and only if $\rho,\sigma$ and $\theta$ satisfy {\rm (C3)}.
\end{pro}

\begin{proof}
If we assume $C(\rho,\sigma;\theta)$, then (C3) is true by Proposition \ref{TCOpsti}. Now assume that $\rho,\sigma,\theta$ are congruences on $\bS$ such that (C3) is true. According to Proposition \ref{TCOpsti}, it suffices to show that the conditions (C1) and (C2) are also satisfied. To prove (C1) we let $a,b,c,d\in S$ be such that $a\;\rho\;b$ and $c\;\sigma\;d$. Assume that $ac\;\theta\;ad$, then from compatibility of $\theta$ we obtain $ba^{-1}ac\;\theta\;ba^{-1}ad$. Recall that $a\;\rho\;b$ implies $a^{-1}\;\rho\;b^{-1}$, hence $a^{-1}a\;\rho\;b^{-1}b$. Now from (C3) we obtain that $b(a^{-1}a)c\;\theta\;b(a^{-1}a)d$ implies $bc=b(b^{-1}b)c\;\theta\;b(b^{-1}b)d=bd$. Therefore, the condition (C1) is satisfied. Analogously, we can prove (C2). By Proposition \ref{TCOpsti}, the condition $C(\rho,\sigma;\theta)$ is true.
\end{proof}

In the next proposition we obtain that similar equivalent conditions for $C(\alpha,\beta;\delta)$ are true for congruences on groups. It follows directly from the fact that every group $\bG$ is a completely simple semigroup, $^{-1}$ is a unary operation compatible with all congruences on the group $\bG$, and Proposition \ref{PropAFu}.

\begin{pro}\label{TCOpstiGrupe}
Let $\mathbf{G}=(G,\cdot,^{-1},e)$ be a group, and let $\alpha,\beta,\delta\in\Con(\mathbf{G})$. Then $C(\alpha,\beta;\delta)$ if and only if $\alpha,\beta$ and $\delta$ satisfy {\rm (C3)}.
\end{pro}

It is useful to observe that in completely simple semigroups represented by normalized Rees matrix semigroups, the condition (C3) can be further simplified. We will say that congruences $\rho,\sigma,\delta$ of $\bS=\mathcal{M}[G;I,\Lambda;P]$ satisfy the condition
\begin{enumerate}
\item[(C3.1)] if $i_1\;\rho_I\;i_2$ then
$$c_1(i_1,e,1)c_2\;\delta\;d_1(i_1,e,1)d_2 \implies c_1(i_2,e,1)c_2\;\delta\;d_1(i_2,e,1)d_2;$$
\item[(C3.2)] if $f_1\;\rho_G\;f_2$ then
$$c_1(1,f_1,1)c_2\;\delta\;d_1(1,f_1,1)d_2 \implies c_1(1,f_2,1)c_2\;\delta\;d_1(1,f_2,1)d_2;$$
\item[(C3.3)] if $\lambda_1\;\rho_{\Lambda}\;\lambda_2$ then
$$c_1(1,e,\lambda_1)c_2\;\delta\;d_1(1,e,\lambda_1)d_2 \implies c_1(1,e,\lambda_2)c_2\;\delta\;d_1(1,e,\lambda_2)d_2;$$
\end{enumerate}
for all $i_1,i_2\in I$, $f_1,f_2\in G$, $\lambda_1,\lambda_2\in\Lambda$ and $c_s,d_s\in S$ such that $c_s\;\sigma\;d_s$, $s=1,2$.

\begin{pro}\label{C3Alt}
Let $\bS=\mathcal{M}[G;I,\Lambda;P]$ be a completely simple semigroup and let $\rho,\sigma,\delta\in\Con(\bS)$. Then the condition $C(\rho,\sigma;\delta)$ is true if and only if $\rho,\sigma$ and $\delta$ satisfy the conditions \rm{(C3.1)}, \rm{(C3.2)} and \rm{(C3.3)}.
\end{pro}

\begin{proof}
By Proposition \ref{TCuslovRMSC3} we know that condition $C(\rho,\sigma;\delta)$ is equivalent to the condition (C3) for congruences $\rho,\sigma,\delta$. Hence, it suffices to show that (C3) is equivalent to the conditions (C3.1), (C3.2) and (C3.3). Note that condition (C3) directly implies conditions (C3.1), (C3.2) and (C3.3), since they can be regarded as special cases of (C3). On the other hand, assume that congruences $\rho,\sigma$ and $\delta$ satisfy the conditions (C3.1), (C3.2) and (C3.3). Let $a,b\in S$, $c_s,d_s\in S$ be such that $a\;\rho\;b$, $c_s\;\sigma\;d_s$, $s=1,2$, where we have denoted $a=(i_1,f_1,\lambda_1)$, $b=(i_2,f_2,\lambda_2)$. Note that we can write $(i_r,f_r,\lambda_r)$ as the product $(i_r,e,1)(1,f_r,1)(1,e,\lambda_r)$, for $r=1,2$. Relation $(i_1,f_1,\lambda_1)\;\rho\;(i_2,f_2,\lambda_2)$ is equivalent to $i_1\;\rho_I\;\i_2$, $f_1\;\rho_G\;f_2$ and $\lambda_1\;\rho_{\Lambda}\;\lambda_2$. From definition of $\rho_I$, $\rho_G$ and $\rho_{\Lambda}$ we obtain $(i_1,e,1)\;\rho\;(i_2,e,1)$, $(1,f_1,1)\;\rho\;(1,f_2,1)$ and $(1,e,\lambda_1)\;\rho\;(1,e,\lambda_2)$, respectively. Now assume that $c_1ac_2\;\delta\;d_1ad_2$, then we have
\begin{equation}
c_1 (i_1,e,1) \Bigl( (1,f_1,1) (1,e,\lambda_1) c_2 \Bigr) \;\delta\;  d_1 (i_1,e,1) \Bigl( (1,f_1,1)(1,e,\lambda_1) d_2 \Bigr). \label{eqn1:C3Alt}
\end{equation}
From $c_2\;\sigma\;d_2$ we obtain $(1,f_1,1) (1,e,\lambda_1) c_2\;\sigma\;(1,f_1,1)(1,e,\lambda_1) d_2$, using compatibility of the congruence $\sigma$. Now by condition (C3.1) it follows that
\begin{equation}
\Bigl( c_1 (i_2,e,1) \Bigr) (1,f_1,1) \Bigl( (1,e,\lambda_1) c_2 \Bigr) \;\delta\; \Bigl( d_1 (i_2,e,1) \Bigr) (1,f_1,1) \Bigl( (1,e,\lambda_1) d_2 \Bigr). \label{eqn2:C3Alt}
\end{equation}
Further on, relations $c_1\;\sigma\;d_1$ and $c_2\;\sigma\;d_2$ imply $c_1(i_2,e,1)\;\sigma\;d_1(i_2,e,1)$ and $(1,e,\lambda_1)c_2\;\sigma\;(1,e,\lambda_1)d_2$, respectively. Applying the condition (C3.2) on \eqref{eqn2:C3Alt} we obtain
\begin{equation}
\Bigl( c_1 (i_2,e,1) (1,f_2,1) \Bigr) (1,e,\lambda_1) c_2  \;\delta\; \Bigl( d_1 (i_2,e,1)(1,f_2,1) \Bigr) (1,e,\lambda_1) d_2. \label{eqn3:C3Alt}
\end{equation}
Relation $c_1\;\sigma\;d_1$ also gives us $c_1(i_2,e,1) (1,f_2,1)\;\sigma\;d_1(i_2,e,1)(1,f_2,1)$. Finally, by the condition (C3.3) relation \eqref{eqn3:C3Alt} implies
\begin{align*}
c_1bc_2 = c_1 (i_2,e,1) (1,f_2,1) (1,e,\lambda_2) c_2 \;\delta\; d_1 (i_2,e,1)(1,f_2,1)(1,e,\lambda_2) d_2 = d_1bd_2.
\end{align*}
Hence, congruences $\rho,\sigma$ and $\delta$ satisfy the condition (C3).
\end{proof}

\section{Commutator on Completely Simple Semigroups}

\begin{lem}\label{LemaLinkedTriple}
Let $\bS=\mathcal{M}[G;I,\Lambda;P]$ be a completely simple semigroup. Let $\rho,\sigma\in\Con(\bS)$ and let $\tau$ be a congruence of $\mathbf{G}$. Then $(0_I,\tau,0_{\Lambda})$ is a linked triple.
\end{lem}

\begin{proof}
Proof follows directly by checking the defining conditions for linked triples (see \cite[page 91]{Howie}).
%Let $i,j\in I$ and $\lambda,\mu\in\Lambda$. First assume that $i\;0_I\;j$, that is $i=j$. Then $q_{\lambda \mu i j}= p_{\lambda i} p_{\mu i}^{-1} p_{\mu j} p_{\lambda j}^{-1}= p_{\lambda i} p_{\mu i}^{-1} p_{\mu i} p_{\lambda i}^{-1} = e$, therefore, the condition $q_{\lambda \mu i j}\;\tau\;e$ is trivially true. Similarly, if $\lambda\;0_{\Lambda}\;\mu$, that is $\lambda=\mu$, then $q_{\lambda \mu i j}=e$, which implies $q_{\lambda \mu i j}\;\tau\;e$. Hence, $(0_I,\tau,0_{\Lambda})$ is a linked triple.
\end{proof}

\begin{pro}\label{PropC(rho,sigma;H)}
Let $\bS=\mathcal{M}[G;I,\Lambda;P]$ be a completely simple semigroup and let $\rho,\sigma$ be congruences on $\bS$. Then the commutator $[\rho,\sigma]$ is determined by the linked triple $(0_I,[\rho,\sigma]_G,0_{\Lambda})$.
\end{pro}

\begin{proof}
By \cite[Lemma 3.7]{RM23} we obtain $C(\rho,\sigma;\mathcal{H})$, because $\mathcal{H}_I=0_I$ and $\mathcal{H}_{\Lambda}=0_{\Lambda}$. By Definition \ref{McMcDefn4.150} and Proposition \ref{PropCongOrd} we obtain the statement.
\end{proof}

\begin{pro}\label{PropC(rho_G,sigma_G;[]}
Let $\bS=\mathcal{M}[G;I,\Lambda;P]$ be a completely simple semigroup. If $\rho,\sigma$ are congruences on $\bS$, then $[\rho_G,\sigma_G]\leq[\rho,\sigma]_G$.
\end{pro}

\begin{proof}
By Proposition \ref{TCOpstiGrupe}, it suffices to prove that the condition (C3) is satisfied for congruences $\rho_G,\sigma_G$ and $[\rho,\sigma]_G$ on the group $\mathbf{G}$. Let $a,b,c_s,d_s\in G$, $s=1,2$ be such that $a\;\rho_G\;b$, $c_s\;\sigma_G\;d_s$, $s=1,2$, and $c_1ac_2\;[\rho,\sigma]_G\;d_1ad_2$. From the definition of congruences $\rho_G,\sigma_G$ and $[\rho,\sigma]_G$ it follows that these relations are equivalent to $(1,a,1)\;\rho\;(1,b,1)$, $(1,c_s,1)\;\sigma\;(1,d_s,1)$, $s=1,2$, and $(1,c_1ac_2,1)\;[\rho,\sigma]\;(1,d_1ad_2,1)$. Using the identity $p_{11}=e$, we can write the last relation as $(1,c_1,1)(1,a,1)(1,c_2,1) \;[\rho,\sigma]\; (1,d_1,1)(1,a,1)(1,d_2,1)$. Now, from the condition $C(\rho,\sigma;[\rho,\sigma])$ we obtain
\begin{equation*}
(1,c_1bc_2,1) = (1,c_1,1)(1,b,1)(1,c_2,1) \;[\rho,\sigma]\; (1,d_1,1)(1,b,1)(1,d_2,1) = (1,d_1bd_2,1),
\end{equation*}
that is, $c_1bc_2\;[\rho,\sigma]_G\;d_1bd_2$, which proves the condition (C3), and consequently $C(\rho_G,\sigma_G;[\rho,\sigma]_G)$. Finally, the desired inequality follows from the Definition \ref{McMcDefn4.150}.
\end{proof}

\begin{de}\label{thetarhosigma}
Let $\bS=\mathcal{M}[G;I,\Lambda;P]$ be a completely simple semigroup and let $\rho,\sigma\in\Con(\bS)$. By $\Theta_{\rho,\sigma}$ we denote the congruence on $\bG$ generated by all ordered pairs $(p_{\mu i} p_{\lambda i}^{-1},p_{\mu j}p_{\lambda j}^{-1})$ where $i,j\in I$ and $\lambda,\mu\in\Lambda$ are such that $i\;\rho_I\;j$ and $\lambda\;\sigma_{\Lambda}\;\mu$, or $\lambda\;\rho_{\Lambda}\;\mu$ and $i\;\sigma_I\;j$.
\end{de}

\begin{rem}\label{RemTheta}
Note that if $\rho=1_S=\sigma$, then the congruence $\Theta_{1,1}$ is generated by all pairs $(p_{\mu i}p_{\lambda i}^{-1},p_{\mu j}p_{\lambda j}^{-1})$, where $i,j\in I$, $\lambda,\mu\in\Lambda$. Since $p_{\lambda 1}=e=p_{1 1}=p_{1 i}$, it follows that $p_{\lambda i}^{-1}\Theta_{1,1}\;e$ for every $i\in I$, $\lambda\in\Lambda$. Therefore, $\Theta_{1,1}$ is the congruence on $\bG$ generated by the set $\mathcal{P}=\{p_{\lambda i}:\lambda\in\Lambda,i\in I\}$. In particular, if $\mathcal{P}$ contains a complete set of generators of the group $\bG$, then $\Theta_{1,1}=1_G$.
\end{rem}

\begin{pro}\label{CorTheta}
Let $\bS=\mathcal{M}[G;I,\Lambda;P]$ be a completely simple semigroup and let $\rho,\sigma$ be congruences on $\bS$. Then $\Theta_{\rho,\sigma}\vee[\rho_G,\sigma_G] \leq [\rho,\sigma]_G$.
\end{pro}

\begin{proof}
Since $C(\rho,\sigma;[\rho,\sigma])$, by Proposition \ref{TCOpsti} it follows that the conditions (C1) and (C2) are true for congruences $\rho,\sigma$ and $[\rho,\sigma]$. Let $i,j\in I$ and $\lambda,\mu\in\Lambda$ be such that $\lambda\;\rho_{\Lambda}\;\mu$ and $i\;\sigma_I\;j$, or $i\;\rho_I\;j$ and $\lambda\;\sigma_{\Lambda}\;\mu$. In the first case, we obtain $a:=(i,p_{\lambda i}^{-1},\lambda)\;\rho\;(i,p_{\mu i}^{-1},\mu)=:b$ and $(i,p_{\lambda i}^{-1},\lambda)\;\sigma\;(j,p_{\lambda j}^{-1},\lambda)=:d$, by definition of $\rho_{\Lambda}$ and $\sigma_I$ respectively. If we denote $c:=a$, then we can easily check that $ac=a^2=a=ad$. Hence, we obtain $ac\;[\rho,\sigma]\;ad$, so the condition (C1) implies $bc\;[\rho,\sigma]\;bd$, that is,
\begin{align*}
(i,p_{\lambda i}^{-1},\lambda) = (i,p_{\mu i}^{-1},\mu)(i,p_{\lambda i}^{-1},\lambda) \;[\rho,\sigma]\; (i,p_{\mu i}^{-1},\mu)(j,p_{\lambda j}^{-1},\lambda)= (i,p_{\mu i}^{-1}p_{\mu j} p_{\lambda j}^{-1},\lambda).
\end{align*}
From previous relation we have $p_{\lambda i}^{-1}\;[\rho,\sigma]_G\; p_{\mu i}^{-1} p_{\mu j} p_{\lambda j}^{-1}$, that is $p_{\mu i} p_{\lambda i}^{-1}\;[\rho,\sigma]_G\; p_{\mu j} p_{\lambda j}^{-1}$. Dually, in the second case we also obtain $p_{\mu i} p_{\lambda i}^{-1}\;[\rho,\sigma]_G\; p_{\mu j} p_{\lambda j}^{-1}$. Hence, by Definition \ref{thetarhosigma} we have the inequality $\Theta_{\rho,\sigma}\leq [\rho,\sigma]_G$. From  Proposition \ref{PropC(rho_G,sigma_G;[]} it follows that $\Theta_{\rho,\sigma}\vee[\rho_G,\sigma_G] \leq [\rho,\sigma]_G$.
\end{proof}

\begin{lem}\label{LemmaComm}
Let $\bS=\mathcal{M}[G;I,\Lambda;P]$ be a completely simple semigroup and let $\rho,\sigma,\delta\in\Con(\bS)$. Then $C(\rho,\sigma;\delta)$ if and only if $\Theta_{\rho,\sigma}\vee[\rho_G,\sigma_G] \leq\delta_G$.
\end{lem}

\begin{proof}
($\rightarrow$) From Proposition \ref{CorTheta} it follows that $\Theta_{\rho,\sigma}\vee[\rho_G,\sigma_G] \leq [\rho,\sigma]_G$. By the definition of the commutator, condition $C(\rho,\sigma;\delta)$ implies $[\rho,\sigma]\leq\delta$. Proposition \ref{PropCongOrd} then gives us $[\rho,\sigma]_G\leq\delta_G$, which further implies $\Theta_{\rho,\sigma}\vee[\rho_G,\sigma_G]\leq \delta_G$.

($\leftarrow$) Assume now that $\Theta_{\rho,\sigma}\vee[\rho_G,\sigma_G]\leq \delta_G$. Since groups form a congruence modular variety, by Proposition \ref{PropCongModVar} the inequality $[\rho_G,\sigma_G]\leq\delta_G$ implies $C(\rho_G,\sigma_G;\delta_G)$. By Proposition \ref{C3Alt} it suffices to prove that conditions (C3.1), (C3.2) and (C3.3) are true for $\rho,\sigma$ and $\delta$. Let $i_1,i_2\in I$, $f_1,f_2\in G$, $\lambda_1,\lambda_2\in\Lambda$ and let $c_s,d_s \in S$, $s=1,2$, be such that $c_s\;\sigma\;d_s$, $s=1,2$. We will use the notation $c_s=(j_s,g_s,\mu_s)$, $d_s=(k_s,h_s,\nu_s)$, for $s=1,2$. Note that from $c_1\;\sigma\;d_1$ and $c_2\;\sigma\;d_2$ we obtain, respectively,
\begin{equation}
g_1\;\sigma_G\;h_1 \mbox{ and } g_2\;\sigma_G\;h_2. \label{eqn1:Lemma4.7}
\end{equation}

\noindent (C3.1) Assume that $i_1\;\rho_I\;i_2$ and $c_1(i_1,e,1)c_2\;\delta\;d_1(i_1,e,1)d_2$, that is
\begin{equation}
(j_1,g_1p_{\mu_1 i_1}g_2,\mu_2) \;\delta\; (k_1,h_1p_{\nu_1 i_1}h_2,\nu_2). \label{eqn2:Lemma4.7}
\end{equation}
From definition of $\rho_I$ we obtain that $(i_1,p_{\mu_1 i_1}^{-1},\mu_1)\;\rho\;(i_2,p_{\mu_1 i_2}^{-1},\mu_1)$. This implies $p_{\mu_1 i_1}^{-1}\;\rho_G\;p_{\mu_1 i_2}^{-1}$, that is
\begin{equation}
p_{\mu_1 i_1} \;\rho_G\; p_{\mu_1 i_2}. \label{eqn3:Lemma4.7}
\end{equation}
By multiplying relation $(j_1,g_1,\mu_1)=c_1\;\sigma\;d_1=(k_1,h_1,\nu_1)$ on the right with $(i_1,e,1)$ we obtain $(j_1,g_1p_{\mu_1 i_1},1)\;\sigma\;(k_1,h_1p_{\nu_1 i_1},1)$. From previous relation we have $g_1p_{\mu_1 i_1}\;\sigma_G\;h_1p_{\nu_1 i_1}$, that gives us
\begin{equation}
g_1\;\sigma_G\;h_1p_{\nu_1 i_1} p_{\mu_1 i_1}^{-1}. \label{eqn4:Lemma4.7}
\end{equation}
On the other hand, from relation $c_1\;\sigma\;d_1$ we obtain $\mu_1\;\sigma_{\Lambda}\;\nu_1$. We have $i_1\;\rho_I\;i_2$, hence by Definition \ref{thetarhosigma} it follows that $p_{\nu_1 i_1} p_{\mu_1 i_1}^{-1} \;\Theta_{\rho,\sigma}\; p_{\nu_1 i_2} p_{\mu_1 i_2}^{-1}$, and therefore $p_{\nu_1 i_1}p_{\mu_1 i_1}^{-1} p_{\mu_1 i_2}\;\Theta_{\rho,\sigma}\;p_{\nu_1 i_2}$. Since $\Theta_{\rho,\sigma}\leq \delta_G$, we obtain $p_{\nu_1 i_1}p_{\mu_1 i_1}^{-1} p_{\mu_1 i_2}\;\delta_G\;p_{\nu_1 i_2}$, which further implies
\begin{equation}
h_1p_{\nu_1 i_1}p_{\mu_1 i_1}^{-1} p_{\mu_1 i_2}h_2\;\delta_G\;h_1p_{\nu_1 i_2}h_2. \label{eqn5:Lemma4.7}
\end{equation}
From \eqref{eqn2:Lemma4.7} we have $j_1\;\delta_I\;k_1$, $\mu_2\;\delta_{\Lambda}\;\nu_2$ and $g_1p_{\mu_1 i_1}g_2\;\delta_G\;h_1p_{\nu_1 i_1}h_2$. The last relation is equivalent to
\begin{equation}
g_1p_{\mu_1 i_1}g_2\;\delta_G\; (h_1p_{\nu_1 i_1} p_{\mu_1 i_1}^{-1}) p_{\mu_1 i_1} h_2. \label{eqn6:Lemma4.7}
\end{equation}
From \eqref{eqn4:Lemma4.7}, \eqref{eqn3:Lemma4.7}, \eqref{eqn1:Lemma4.7},  and \eqref{eqn6:Lemma4.7} we obtain
\begin{equation}
g_1 p_{\mu_1 i_2}  g_2 \;\delta_G\; (h_1 p_{\nu_1 i_1} p_{\mu_1 i_1}^{-1}) p_{\mu_1 i_2} h_2, \label{eqn7:Lemma4.7}
\end{equation}
by condition (C3) for $\rho_G,\sigma_G$ and $\delta_G$, because of $C(\rho_G,\sigma_G;\delta_G)$ and Proposition \ref{TCOpstiGrupe}. Using transitivity of $\delta_G$ on \eqref{eqn7:Lemma4.7} and \eqref{eqn5:Lemma4.7} we have $g_1p_{\mu_1 i_2}g_2\;\delta_G\;h_1p_{\nu_1 i_2}h_2$. Since $j_1\;\delta_I\;k_1$, $\mu_2\;\delta_{\Lambda}\;\nu_2$, it follows that
\begin{equation*}
c_1(i_2,e,1)c_2 = (j_1,g_1p_{\mu_1 i_2} g_2,\mu_2) \;\delta\; (k_1, h_1 p_{\nu_1 i_2} h_2, \nu_2) = d_1(i_2,e,1)d_2.
\end{equation*}

\noindent (C3.2) Assume that $f_1\;\rho_G\;f_2$ and $c_1(1,f_1,1)c_2\;\delta\;d_1(1,f_1,1)d_2$, that can be written as $(j_1,g_1f_1g_2,\mu_2)\;\delta\;(k_1,h_1f_1h_2,\nu_2)$.
The last relation is equivalent to $j_1\;\delta_I\;k_1$, $\mu_2\;\delta_{\Lambda}\;\nu_2$ and $g_1f_1g_2\;\delta_G\;h_1f_1h_2$. From $f_1\;\rho_G\;f_2$, \eqref{eqn1:Lemma4.7} and  $g_1f_1g_2\;\delta_G\;h_1f_1h_2$ we obtain $g_1f_2g_2\;\delta_G\;h_1f_2h_2$, by condition (C3). Using relations $j_1\;\delta_I\;k_1$, $\mu_2\;\delta_{\Lambda}\;\nu_2$, we obtain $c_1bc_2=(j_1,g_1f_2g_2,\mu_2)\;\delta\;(k_1,h_1f_2h_2,\nu_2)=d_1bd_2$.

\noindent (C3.3) Dually to (C3.1).
\end{proof}

Now we give the proof of the main theorem.

\noindent \textit{Proof of Theorem \ref{ThmComm}}\\
Denote by $\theta$ the congruence determined by the linked triple $(0_I,\Theta_{\rho,\sigma}\vee[\rho_G,\sigma_G],0_{\Lambda})$. Lemma \ref{LemmaComm} directly implies $C(\rho,\sigma;\theta)$. By the definition of the commutator, it follows that $[\rho,\sigma]\leq\theta$. In Proposition \ref{PropC(rho,sigma;H)} we have proved that the commutator $[\rho,\sigma]$ is determined by the linked triple $(0_I,[\rho,\sigma]_G,0_{\Lambda})$. Proposition \ref{PropCongOrd} then implies that $[\rho,\sigma]_G\leq \theta_G=\Theta_{\rho,\sigma}\vee[\rho_G,\sigma_G]$. On the other hand, Proposition \ref{CorTheta} gives us $\Theta_{\rho,\sigma}\vee[\rho_G,\sigma_G] \leq [\rho,\sigma]_G$. Therefore, $[\rho,\sigma]_G$ is equal to $[\rho_G,\sigma_G]\vee\Theta_{\rho,\sigma}$, and consequently, $[\rho,\sigma] = \theta$ is determined by the linked triple $(0_I,[\rho_G,\sigma_G]\vee\Theta_{\rho,\sigma},0_{\Lambda})$. \qed

If we assume some additional conditions on congruences $\rho$ and $\sigma$, the group part of the commutator $[\rho,\sigma]$ simplifies, as we can see in the following statement. Its significance is justified in the next section.

\begin{cor}\label{CorOfThm1}
Let $\bS=\mathcal{M}[G;I,\Lambda;P]$ be a completely simple semigroup and let $\rho,\sigma$ be congruences on $\bS$. If at least one of the following conditions:
\begin{enumerate}
\item[(i)] $\rho_I=\sigma_I=0_I$;
\item[(ii)] $\rho_{\Lambda} = \sigma_{\Lambda} = 0_{\Lambda}$;
\item[(iii)] $\rho_I=0_I$ and $\rho_{\Lambda}=0_{\Lambda}$;
\item[(iv)] $\sigma_I=0_I$ and $\sigma_{\Lambda}=0_{\Lambda}$.
\end{enumerate}
is satisfied, then $[\rho_G,\sigma_G] = [\rho,\sigma]_G$.
\end{cor}
\begin{proof}
From Theorem \ref{ThmComm} we obtain $[\rho,\sigma]_G = [\rho_G,\sigma_G]\vee\Theta_{\rho,\sigma}$, where $\Theta_{\rho,\sigma}$ is the congruence generated by ordered pairs $(p_{\mu i}p_{\lambda i}^{-1},p_{\mu j}p_{\lambda j}^{-1})$, where $i,j\in I$, $\lambda,\mu\in\Lambda$ are such that $i\;\rho_I\;j$ and $\lambda\;\sigma_{\Lambda}\;\mu$, or $i\;\sigma_I\;j$ and $\lambda\;\rho_{\Lambda}\;\mu$. Note that if either of the conditions: (i)--(iv) is satisfied, then $p_{\mu i}p_{\lambda i}^{-1} = p_{\mu j}p_{\lambda j}^{-1}$ for every $i,j\in I$ and $\lambda,\mu\in\Lambda$ such that $i\;\rho_I\;j$ and $\lambda\;\sigma_{\Lambda}\;\mu$, or $\lambda\;\rho_{\Lambda}\;\mu$ and $i\;\sigma_I\;j$. Hence, in either case the congruence $\Theta_{\rho,\sigma}$ is equal to $0_G$, which further implies $[\rho,\sigma]_G=[\rho_G,\sigma_G]$.
\end{proof}

\section{Applications and Examples}

\begin{exa}\label{Example}
We will determine the commutator $[\rho,\sigma]$ when $\rho,\sigma$ are congruences $\Hg,\Lg$ and $\Rg$ on a completely simple semigroups $\bS=\mathcal{M}[G;I,\Lambda;P]$. \\
(a) If $\rho=\Hg=\sigma$, then we can apply the Corollary \ref{CorOfThm1}(iii), which gives us $[\Hg,\Hg]_G = [\Hg_G,\Hg_G] = [1_G,1_G]$. Therefore, the commutator of $[\Hg,\Hg]$ is equal to $(0_I,[1_G,1_G],0_{\Lambda})$. \\
(b) In the case $\rho=\sigma=\Lg$ we can apply the Corollary \ref{CorOfThm1}(ii), which gives us $[\Lg,\Lg]_G=[\Lg_G,\Lg_G]=[1_G,1_G]$. Hence, the commutator $[\Lg,\Lg]$ corresponds to $(0_I,[1_G,1_G],0_{\Lambda})$. \\
(c) Similarly to (b), from Corollary \ref{CorOfThm1}(i) it follows that $[\Rg,\Rg]$ also corresponds to $(0_I,[1_G,1_G],0_{\Lambda})$.\\
(d) By Theorem \ref{ThmComm} we have $[\Lg,\Rg]_G = [\Lg_G,\Rg_G]\vee\Theta_{\Lg,\Rg} = [1_G,1_G]\vee\Theta_{\Lg,\Rg}$. Since $\Rg_I=0_I$, $\Rg_{\Lambda}=1_{\Lambda}$ and $\Lg_I=1_I$, $\Lg_{\Lambda}=0_{\Lambda}$, from Definition \ref{thetarhosigma} it fo\-llows that $\Theta_{\Lg,\Rg}$ is the congruence generated by all ordered pairs $(p_{\mu i}p_{\lambda i}^{-1},p_{\mu j}p_{\lambda j}^{-1})$, where $i,j\in I$, $\lambda,\mu\in\Lambda$ are arbitrary. Hence, $\Theta_{\Lg,\Rg}=\Theta_{1,1}$, and consequently $[\Lg,\Rg]_G=[1_G,1_G]\vee\Theta_{1,1}=[1_S,1_S]_G$, that is $[\mathcal{L},\mathcal{R}]=[1_S,1_S]$.
\end{exa}

Let us recall that congruence lattice of a completely simple semigroup is semi-modular \cite[Theorem 3.6.2]{Howie}. However, Lemma \ref{LemmaComm} and Theorem \ref{ThmComm} allows us to prove some properties of the commutator which are satisfied in modular varieties.

\begin{pro}\label{Cor1}
Let $\bS$ be a completely simple semigroup. Let $\rho,\sigma,\delta\in\Con(\bS)$, then $C(\rho,\sigma;\delta)$ if and only if $[\rho,\sigma]\leq \delta$.
\end{pro}

\begin{proof}
Follows directly from Lemma \ref{LemmaComm} and Theorem \ref{ThmComm}.
\end{proof}

\begin{pro}\label{PropFactor}
Let $\bS$ be a completely simple semigroup, and let $\rho,\sigma,\eta\in\Con(\bS)$ be such that $\eta\leq\rho,\sigma$. Then in $S/\eta$ we have $[\rho/\eta,\sigma/\eta] = ([\rho,\sigma]\vee\eta)/\eta$.
\end{pro}

\begin{proof}
From Definition \ref{McMcDefn4.148} one can see that in arbitrary algebra for all congruences $\rho,\sigma,\theta$ and $\eta$ such that $\eta\leq\rho,\sigma,\theta$, we have $C(\rho,\sigma;\theta)$ if and only if $C(\rho/\eta,\sigma/\eta;\theta/\eta)$. Hence, whenever $C(\rho/\eta,\sigma/\eta;\theta/\eta)$ is true for some $\theta/\eta\in\Con(\bS/\eta)$, then $[\rho,\sigma]\vee\eta\leq\theta$ and therefore $([\rho,\sigma]\vee\eta)/\eta \leq \theta/\eta$. By Proposition \ref{Cor1}, from $[\rho,\sigma]\leq [\rho,\sigma]\vee\eta$ we obtain $C(\rho,\sigma;[\rho,\sigma]\vee\eta)$, that is equivalent to $C(\rho/\eta,\sigma/\eta; ([\rho,\sigma]\vee\eta)/\eta)$, and hence $[\rho/\eta,\sigma/\eta] = ([\rho,\sigma]\vee\eta)/\eta$, by Definition \ref{McMcDefn4.150}.
\end{proof}

\begin{pro}\label{Cor[a,b]=[b,a]}
Let $\bS$ be a completely simple semigroup. Let $\rho,\sigma\in\Con(\bS)$, then $[\rho,\sigma]=[\sigma,\rho]$.
\end{pro}

\begin{proof}
Let $\mathcal{M}[G;I,\Lambda;P]$ be a Rees matrix representation of $\bS$. We observe that the definition of congruence $\Theta_{\rho,\sigma}$ is symmetric, hence $\Theta_{\rho,\sigma}=\Theta_{\sigma,\rho}$. Since $[\rho_G,\sigma_G]=[\sigma_G,\rho_G]$ in the group $\bG$, by Theorem \ref{ThmComm} we obtain $[\rho,\sigma] = [\sigma,\rho]$.
\end{proof}

\begin{pro}\label{Cor[a,bigveeb]=bigvee[a,b]}
Let $\bS$ be a completely simple semigroup. Let $\rho$ and $\sigma_i$, $i\in I$ be congruences on $\bS$. Then $[\rho,\bigvee_{i\in I} \sigma_i] = \bigvee_{i\in I} [\rho,\sigma]$.
\end{pro}

\begin{proof}
The proof of this proposition follows the same steps as the proof of \cite[ Proposition 4.3]{FM:CTCMV}. From Proposition \ref{CommutatorProperties} we obtain $\bigvee_{i\in I} [\rho,\sigma_i] \leq [\rho,\bigvee_{i\in I}\sigma_i]$. In order to prove $C(\rho,\bigvee_{i\in I}\sigma_i;\bigvee_{i\in I}[\rho,\sigma_i])$ we use the Proposition \ref{Cor1} in place of \cite[Proposition 4.2]{FM:CTCMV}. Therefore $[\rho,\bigvee_{i\in I}\sigma_i] \leq \bigvee_{i\in I} [\rho,\sigma_i]$, which completes the proof of the equality from the statement.
\end{proof}

\begin{de}\textup{(cf.\cite{FM:CTCMV}, Definition 6.1)}\label{DefNilp}
Let $\bS$ be a semigroup, and let $\alpha,\beta\in\Con(\bS)$. We define the series of congruences $(\alpha,\beta]^{(k)}=\underbrace{[\alpha,[\alpha,\dotsc,[}_{k}\alpha,\beta]\dotsc]]$, for $k\in\N$. Semigroup $\bS$ is \textit{$n$-nilpotent}, $n\in\N$ if we have the equality $(1_S,1_S]^{(n)}=0_S$. The smallest such $n$ is called the \emph{degree of nilpotency} of the semigroup $\bS$. \\
Similarly, we define a series of congruences $[\alpha]^{(k)}$, $k\in\N$ where $[\alpha]^{(1)} = [\alpha,\alpha]$ and $[\alpha]^{(k)} = [[\alpha,\alpha]^{(k-1)},[\alpha,\alpha]^{(k-1)}]$ for  $k\geq 2$. Semigroup $\bS$ is \textit{$n$-solvable}, $n\in\N$ if we have the equality $[1_S]^{(n)} = 0_S$. The smallest $n$ such that $\bS$ is $n$-solvable is called the \emph{degree of solvability} of the semigroup $\bS$.
\end{de}

\begin{rem}
Notion of nilpotency from Definition \ref{DefNilp} is often refereed to as \emph{left nilpotency}, since we can dually define the series of congruences $[\alpha,\beta)^{(k)}$, $k\in\N$ and a notion of \emph{right nilpotency}. Degrees of left and right nilpotency are independent in the general case  (\cite{KK}). However, from Proposition \ref{Cor[a,b]=[b,a]}, inductively it follows that congruence series $(\rho,\sigma]^{(n)}$ and $[\rho,\sigma)^{(n)}$ are equal, therefore, in completely simple semigroups, notions of left $n$-nilpotency and right $n$-nilpotency are equivalent.
\end{rem}

\begin{lem}\label{LemmaNilpSolv}
Let $\bS=\mathcal{M}[G;I,\Lambda;P]$ be a completely simple semigroup, then for every $k\in\N$, $k\geq 2$ we have
\begin{enumerate}
\item[(a)] $(1_S,1_S]^{(k)}_G = (1_G,[1_G,1_G]\vee\Theta_{1,1}]^{(k-1)}$;
\item[(b)] $[1_S]^{(k)}_G = [[1_G,1_G]\vee\Theta_{1,1}]^{(k-1)}$.
\end{enumerate}
\end{lem}

\begin{proof}
We will give a proof of (a) and (b) by induction on $k\in\N$, $k\geq 2$. First note that from Theorem \ref{ThmComm} it follows that $[1_S,1_S]$ is determined by the linked triple $(0_I,[1_G,1_G]\vee\Theta_{1,1},0_{\Lambda})$. \\
(a) Since $[1_S,1_S]_I=0_I$ and $[1_S,1_S]_{\Lambda}=0_{\Lambda}$ from Corollary \ref{CorOfThm1}(iv) we obtain $(1_S,1_S]^{(2)}_G=[1_S,[1_S,1_S]]_G = [1_G,[1_S,1_S]_G] = (1_G,[1_G,1_G]\vee\Theta_{1,1}]^{(1)}$. Inductively, by using the equalities $(1_S,1_S]^{(k)}_I=0_I$ and $(1_S,1_S]^{(k)}_{\Lambda}=0_{\Lambda}$ that follow from Theorem \ref{ThmComm}, Corollary \ref{CorOfThm1}(iv) gives us $(1_S,1_S]^{(k)}_G = (1_G,[1_G,1_G]\vee\Theta_{1,1}]^{(k-1)}$.\\
(b) Again using $[1_S,1_S]_I=0_I$ and $[1_S,1_S]_{\Lambda}=0_{\Lambda}$, from Corollary \ref{CorOfThm1}(iii) it follows that $[1_S]^{(2)}_G=[[1_S,1_S],[1_S,1_S]]_G = [[1_S,1_S]_G,[1_S,1_S]_G] = [[1_S,1_S]_G]^{(1)}$. Inductively, using the equalities $[1_S]^{(k)}_I=0_I$ and $[1_S]^{(k)}_{\Lambda}=0_{\Lambda}$ obtained from Theorem \ref{ThmComm}, Corollary \ref{CorOfThm1}(iii) gives us $[1_S]^{(k)}_G = [[1_G,1_G]\vee\Theta_{1,1}]^{(k-1)}$.
\end{proof}

\begin{pro}\label{PropNilpSolvEq}
Let $\bS=\mathcal{M}[G;I,\Lambda;P]$ be a completely simple semigroup. Then for every $n\in\N$, $n\geq 2$ we have
\begin{enumerate}
\item[(a)] if the semigroup $\bS$ is $n$-nilpotent, then the group $\bG$ is also $n$-nilpotent;
\item[(b)] if the semigroup $\bS$ is $n$-solvable, then the group $\bG$ is also $n$-solvable.
\end{enumerate}
\end{pro}

\begin{proof}
(a) Note that for all congruences $\alpha,\beta\in\Con(\bG)$, such that $\alpha\leq\beta$ we have $(1,\alpha]^{(k)}\leq(1,\beta]^{(k)}$ by Proposition \ref{CommutatorProperties}. Then for $\alpha=[1_G,1_G]$ and $\beta=[1_G,1_G]\vee\Theta_{1,1}$, it follows $(1_G,1_G]^{(k)}=(1_G,[1_G,1_G]]^{(k-1)} \leq (1_G,[1_G,1_G]\vee\Theta_{1,1}]^{(k-1)}$. Therefore, from Lemma \ref{LemmaNilpSolv} (a) we obtain $(1_G,1_G]^{(k)} \leq (1_S,1_S]^{(k)}_G$ for every $k\in\N$, $k\geq 2$. Consequently, if we have $(1_S,1_S]^{(n)}=0_S$ for some $n\in\N$, $n\geq 2$ then from Proposition \ref{PropC(rho,sigma;H)} it follows that $(1_S,1_S]^{(n)}_G=0_G$. Hence, $(1_G,1_G]^{(n)}=0_G$, that is, the group $\mathbf{G}$ is $n$-nilpotent.

(b) Similarly as in (a), note that for all $\alpha,\beta\in\Con(\bG)$ such that $\alpha\leq\beta$ we have $[\alpha]^{(k)}\leq[\beta]^{(k)}$ by Proposition \ref{CommutatorProperties}. Therefore, for $\alpha=[1_G,1_G]$ and $\beta=[1_G,1_G]\vee\Theta_{1,1}$, we obtain $[1_G]^{(k)}=[[1_G,1_G]]^{(k-1)} \leq [[1_G,1_G]\vee\Theta_{1,1}]^{(k)}$. Then, from Lemma \ref{LemmaNilpSolv} (b) we obtain $[1_G]^{(k)} \leq [1_S]^{(k)}_G$ for every $k\in\N$, $k\geq 2$. Therefore, if we have $[1_S]^{(n)}=0_S$ for some $n\in\N$, $n\geq 2$, then from Proposition \ref{PropC(rho,sigma;H)} it follows that $[1_S]^{(n)}_G=0_G$. Hence, $[1_G]^{(n)}=0_G$, that is, the group $\mathbf{G}$ is $n$-solvable.
\end{proof}

\begin{pro}\label{PropNilpSolvIneq}
Let $\bS=\mathcal{M}[G;I,\Lambda;P]$ be a completely simple semigroup. Then for every $n\in\N$, $n\geq 2$ we have
\begin{enumerate}
\item[(a)] if the group $\bG$ is $n$-nilpotent, then the semigroup $\bS$ is $(n+1)$-nilpotent;
\item[(b)] if the group $\bG$ is $n$-solvable, then the semigroup $\bS$ is $(n+1)$-solvable.
\end{enumerate}
\end{pro}

\begin{proof}
(a) From Lemma \ref{LemmaNilpSolv}(a) we have $(1_S,1_S]^{(k)}_G = (1_G,[1_G,1_G]\vee\Theta_{1,1}]^{(k-1)}$, for every $k\in\N$, $k\geq 2$. Since $[1_G,1_G]\vee\Theta_{1,1} \leq 1_G$, applying the Proposition \ref{CommutatorProperties} inductively, we obtain that $(1_G,[1_G,1_G]\vee\Theta_{1,1}]^{(k-1)} \leq (1_G,1_G]^{(k-1)}$. Therefore, we have proved
\begin{equation}
(1_S,1_S]^{(k)}_G \leq  (1_G,1_G]^{(k-1)}, \mbox{ for every } k\in\N, k\geq 2. \label{eqn:AA}
\end{equation}

Now assume that the group $\mathbf{G}$ is $n$-nilpotent, where $n\in\N$, $n\geq 2$. Then we have $(1_G,1_G]^{(n)}=0_G$. From inequality \eqref{eqn:AA} it follows that $(1_S,1_S]^{(n+1)}_G \leq 0_G$. Therefore, from Proposition \ref{PropC(rho,sigma;H)} we obtain that $(1_S,1_S]^{(n+1)} = 0_S$, that is, semigroup $\bS$ is $(n+1)$-nilpotent.

(b) By Lemma \ref{LemmaNilpSolv}(b) we have the equality $[1_S]^{(k)}_G = [[1_G,1_G]\vee\Theta_{1,1}]^{(k-1)}$, for every $k\in\N$, $k\geq 2$. From inequality $[1_G,1_G]\vee\Theta_{1,1}\leq 1_G$, by Proposition \ref{CommutatorProperties} we obtain that $[[1_G,1_G]\vee\Theta_{1,1}]^{(k-1)} \leq [1_G]^{(k-1)}$ for every $k\in\N$, $k\geq 2$. Therefore, we have the inequality
\begin{equation}
[1_S]^{(k)}_G \leq  [1_G]^{(k-1)}, \mbox{ for every } k\in\N, k\geq 2. \label{eqn:BB}
\end{equation}

Now assume that the group $\mathbf{G}$ is $n$-solvable, where $n\in\N$, $n\geq 2$. Then we have $[1_G]^{(n)}=0_G$. From inequality \eqref{eqn:BB} it follows that $[1_S]^{(n+1)}_G \leq 0_G$. Therefore, from Proposition \ref{PropC(rho,sigma;H)} we obtain that $[1_S]^{(n+1)} = 0_S$, that is, semigroup $\bS$ is $(n+1)$-solvable.
\end{proof}

\begin{pro}\label{PropNilpSolv}
Let $\bS$ be a completely simple semigroup then:
\begin{enumerate}
\item[(a)] semigroup $\bS$ is nilpotent if and only if its maximal subgroup is a nilpotent group;
\item[(b)] semigroup $\bS$ is solvable if and only if its maximal subgroup is a solvable group.
\end{enumerate}
\end{pro}

\begin{proof}
Proposition \ref{PropNilpSolv} follows directly from Proposition \ref{PropNilpSolvEq} and Proposition \ref{PropNilpSolvIneq}.
\end{proof}

Let us remark that the class of nilpotency (solvability) of the semigroup $\bS=\mathcal{M}[G;I,\Lambda;P]$ and the group $\bG$ can be the same, under a special condition.

\begin{pro}\label{PropNilpSolvCond}
Let $\bS=\mathcal{M}[G;I,\Lambda;P]$ be a completely simple semigroup. If $\Theta_{1,1}\leq [1_G,1_G]$ then for every $n\in\N$, $n\geq 2$ we have
\begin{enumerate}
\item[(a)] if the group $\bG$ is $n$-nilpotent, then the semigroup $\bS$ is $n$-nilpotent;
\item[(b)] if the group $\bG$ is $n$-solvable, then the semigroup $\bS$ is $n$-solvable.
\end{enumerate}
\end{pro}

\begin{proof}
First note that from $\Theta_{1,1}\leq [1_G,1_G]$, we obtain $[1_G,1_G]\vee\Theta_{1,1}=[1_G,1_G]$.\\
(a)  By Lemma \ref{LemmaNilpSolv}(a) it follows $(1_S,1_S]^{(k)}_G = (1_G,[1_G,1_G]]^{(k-1)} = (1_G,1_G]^{(k)}$ for every $k\in\N$, $k\geq 2$. Hence, if the group $\bG$ is $n$-nilpotent, then $(1_S,1_S]^{(n)}_G = (1_G,1_G]^{(n)}=0_G$, and consequently $\bS$ is also $n$-nilpotent.\\
(b) Similarly to (a), from Lemma \ref{LemmaNilpSolv}(b) it follows that $[1_S]^{(k)}_G = [[1_G,1_G]]^{(k-1)} = [1_G]^{(k)}$ for every $k\in\N$, $k\geq 2$. Hence, if the group $\bG$ is $n$-solvable, then $[1_S]^{(n)}_G = [1_G]^{(n)}=0_G$, and consequently $\bS$ is also $n$-solvable.
\end{proof}

\begin{exa}
Let $\bG=D_3$ be the dihedral group. Let $\bS$ be a completely simple semigroup with Rees matrix representation $\bS=\mathcal{M}[G;I,\Lambda;P]$, where $I=\{1,2,3,4\}=\Lambda$, and $P$ is such that $G\subseteq \{p_{\lambda i}:i\in I,\lambda\in\Lambda\}$. By Remark \ref{RemTheta}, we have $\Theta_{1,1}=1_G$. Hence, Theorem \ref{ThmComm} gives us $[1_S,1_S]_G=1_G$. From Lemma \ref{LemmaNilpSolv}(b) we then obtain $[1_S]^{(3)}_G=[1_G]^{(2)}=0_G$, while $[1_S]^{(2)}_G=[1_G,1_G]\neq 0_G$, because $D_3$ is $2$-solvable. Therefore, the semigroup $\bS$ is $3$-solvable, however it is not $2$-solvable. Hence, if the condition $\Theta_{1,1}\leq [1_G,1_G]$ is not satisfied, the degrees of solvability of $\bS$ and $\bG$ do not have to coincide. Similarly, if in the same Rees matrix representation $\bS=\mathcal{M}[G;I,\Lambda;P]$ we take $\bG$ to be the $2$-nilpotent quaternion group $Q_8$, instead of $D_3$, then from Lemma \ref{LemmaNilpSolv}(a) we obtain that the completely simple semigroup $\bS=\mathcal{M}[Q_8;I,\Lambda;P]$ is $3$-nilpotent, but it is not $2$-nilpotent. Therefore, if the condition $\Theta_{1,1}\leq [1_G,1_G]$ is not true, then the degrees of nilpotency of $\bS$ and $\bG$ also do not have to coincide.
\end{exa}

\noindent\textit{Proof of Theorem \ref{ThmNilpSolv}:} \\
(a) ($\rightarrow$) A regular semigroup $\bS$ is nilpotent if and only if it is a nilpotent completely simple semigroup, by \cite[Proposition 3.5]{RM23}. Hence, if we assume that $\bS$ is nilpotent then it must be simple. Recall that all $\mathcal{H}$-classes of a completely simple semigroup $\bS$ are isomorphic to its maximal subgroup $\mathbf{H}$. From Proposition \ref{PropNilpSolv}(a) it follows that $\mathbf{H}$ is a nilpotent group. ($\leftarrow$) Now assume that $\bS$ is a regular simple semigroup, such that all its $\mathcal{H}$-classes are nilpotent groups. From \cite[Proposition 4.1.1]{Howie} it follows that $\bS$ is a completely regular semigroup. Using the fact that $\bS$ is simple, and \cite[Proposition 4.1.2(3)]{Howie}, we obtain that $\bS$ is a completely simple semigroup. Since every $\mathcal{H}$-class of the completely simple semigroup $\bS$ is its maximal subgroup, from Proposition \ref{PropNilpSolv}(a) we obtain that $\bS$ is a nilpotent completely simple semigroup, and therefore nilpotent regular semigroup.\\
(b) Analogous to (a). \qed

\section{Acknowledgments}

We thank to Karl Auinger for useful comments.

%\section{Statements and Declarations}

%The authors declare that they have no relevant financial or non-financial interests to disclose.

\bibliographystyle{plain}

\end{document}